\documentclass[12pt]{amsart}
\usepackage{amsmath}
\usepackage{amsfonts}
\usepackage{amssymb}
\usepackage{amscd}
\usepackage{epic}
\usepackage{eepic}

\newcommand{\printname}[1]
  {\smash{\makebox[0pt]{pace{-2.0in}\raisebox{8pt}{\tiny #1}}}}
\newtheorem {Theorem}   {Theorem} [section]
\newtheorem {Lemma} [Theorem]{Lemma}
\newtheorem {Proposition}[Theorem]{Proposition}

\newtheorem{Definition}[Theorem]{Definition}
\newtheorem{Remark}[Theorem]{Remark}

\newcommand{\ssminus}{{\smallsetminus}}

\newcommand{\nab}[1][]{\ensuremath{\mathrm{\nabla}{#1}}}

\newcommand{\be}{\begin{equation}}
\newcommand{\ee}{\end{equation}}

\newcommand{\lb}[1]{\label{#1}}

\newcommand{\al}{\alpha}

\newcommand{\we}{\wedge}

\newcommand{\lam}{\lambda}

\newcommand{\reff}{\ref}
\newcommand{\Ref}[1]{(\ref{#1})}
\newcommand{\fr}{\frac}

\newcommand{\ra}{\rightarrow}

\newcommand{\ri}{\mathrm{r}}

\newcommand{\Lap}{\Delta}

\newcommand{\wht}{\widehat}
\newcommand{\sk}{SKR}
\def\t{\tau}

\newcommand{\bee}{quasi-Einstein}
\newcommand{\berp}{quasi-Einstein pair}
\newcommand{\ber}{quasi-Einstein metric}
\newcommand{\berf}{quasi-Einstein function}
\newcommand{\berfs}{quasi-Einstein functions}
\newcommand{\Bers}{Quasi-Einstein metrics}
\newcommand{\bers}{quasi-Einstein metrics}
\newcommand{\sol}{Ricci soliton}

\newcommand{\Mt}{M_\t}

\begin{document}

\title{Conformally K\"ahler base metrics for Einstein warped products}
\author[Gideon Maschler]{GIDEON MASCHLER}
\address{Department of Mathematics and Computer Science, Clark University,
Worcester, MA 01610, U.S.A.}
\email{gmaschler@clarku.edu}

\begin{abstract}
A Riemannian metric $\wht{g}$ with Ricci curvature $\wht{\ri}$ is called nontrivial \bee,
in the sense of Case, Shu and Wei, if  it satisfies $(-a/f)\wht{\nab} df+\wht{\ri}=\lambda \wht{g}$,
for a smooth nonconstant function $f$ and constants $\lambda$ and $a>0$. If $a$ is a positive integer, by
a result of Kim and Kim, such a metric forms a base for certain warped Einstein metrics. On a manifold
$M$ of real dimension at least six, let $(g,\t)$ be a pair consisting of a K\"ahler metric $g$ which is
locally K\"ahler irreducible, and a nonconstant Killing potential $\t$. Suppose the metric $\wht{g}=g/\t^2$
is nontrivial \bee\ on $M\setminus\t^{-1}(0)$, and the associated function $f$ is locally a function of $\t$.
Then $(g,\t)$ is an \sk\ pair, a notion defined by Derdzinski and Maschler. This implies
that $M$ is biholomorphic to an open set in the total space of a $CP^1$ bundle whose base manifold admits a
K\"ahler-Einstein metric. If $M$ is additionally compact, it is a total space of such a bundle or complex
projective space. Also,  the function $f$ is affine in $\t^{-1}$ with nonzero constants.
Conversely, in all even dimensions $n\geq 4$, there exist \sk\ pairs $(g,\t)$ and corresponding nonzero
constants $K$ and $L$ for which $g/\t^2$ is nontrivial quasi-Einstein with $f=K\t^{-1}+L$. Additionally,
a result of Case, Shu and Wei on the K\"ahler reducibility of nontrivial K\"ahler \bers\ is reproduced
in dimension at least six in a more explicit form.
\end{abstract}

\maketitle

\setcounter{Theorem}{0}
\renewcommand{\theTheorem}{\Alph{Theorem}}
\renewcommand{\theequation}{\arabic{section}.\arabic{equation}}

\section{Introduction}
\setcounter{equation}{0}

On a manifold $M$ of dimension $n$, the $a$-Bakry-Emery Ricci tensor of a pair $(g,u)$, consisting
of a Riemannian metric $g$ and a smooth function $u$, is defined to be the symmetric $2$-tensor
\[\text{$\ri_u^a=\ri+\nab d u-a^{-1}du\otimes du$, \quad for a constant $0<a<\infty$,}\]
where $\ri$ is the Ricci tensor of $g$ and $\nab d u$ is the Hessian of $u$. If $g$ satisfies the equation $\ri_u^a=\lambda g$,
for $\lam\in \mathbf{R}$, $g$ was called in \cite{csw} a \ber. [This is in contrast with the
use of the term in the physics literature to denote a gradient \sol. Other usages, especially referring to an equation involving
$g$, $\ri$ and $du\otimes du$ {\em but not  $\nab d u$}, also exist \cite{gv}.] The limiting value $a=\infty$, where
a \ber\ becomes a gradient \sol, will not concern us in this work, and thus, in opposition to the convention of \cite{csw}, we
exclude it from the definition.

An important observation made in \cite{csw} involves the substitution $f=\exp(-u/a)$ (for $a$ finite), which converts
the equation for a \ber\ to the form \be\lb{ber}(-a/f)\nab df+\ri=\lambda g.\end{equation} We will call $(g,f)$ a \berp,
and $f$ the \berf. If $f$ is constant, $g$ is Einstein. If \Ref{ber} is satisfied for
some nonconstant $f$, the \berp\ is said to be nontrivial.

The main application of Equation \Ref{ber} relates to the topic of warped Einstein metrics, addressed by Besse \cite{bes}.
According to \cite{kk} (see also \cite{csw}), $(g,f)$  satisfies \Ref{ber} for a positive {\em integer} $a$ on $M$, if and
only if the warped product $M \times_f F^a$ is Einstein, where $F$ is an $a$-dimensional Einstein manifold with Einstein
constant equal to $f\Delta f+(a-1)|\nab f|^2+\lam f^2$. Thus, existence of \bers,  with $a$ a positive integer, is equivalent
to the existence of certain warped Einstein metrics. On compact manifolds, examples of the latter, and therefore, also of the
former, were obtained by L\"u, Page and Pope \cite{lpp}.

Another result obtained in \cite{csw} states that there are no nontrivial K\"ahler \bers\ (with, of course, finite $a$)
on compact manifolds. This result depends on a structure theorem obtained there in the complete simply connected
case, but the argument given is essentially local. In an appendix we give, under mild assumptions, another, somewhat
more explicit version of this structure theorem. However, our main purpose is to address the issue of whether a \ber\ can
be conformal to a K\"ahler metric. We answer this question in the affirmative. In fact, using methods akin to those in \cite{confsol},
we give the following classification and existence theorem which relies on the notion of an \sk\ pair \cite{confsol},
a well understood metric type first studied in \cite{local} (see Definition \reff{skrp-def}).
\begin{Theorem}\lb{A}
On a manifold $M$ of real dimension $n\geq 6$, suppose $g$ is a K\"ahler, and not a local product of
K\"ahler metrics in any neighborhood of some point. Suppose $\t$ is a nonconstant Killing potential for $g$,
and there exists a nontrivial quasi-Einstein pair of the form $(g/\t^2, f(\t))$
on $M\setminus\t^{-1}(0)$.
Then:
\newcounter{ctr}
\begin{list}{}{\usecounter{ctr}}
\item[] the pair $(g,\t)$ is an \sk\ pair;
\item[] the manifold $M$ is biholomorphic to an open set in the total space of a $CP^1$ bundle whose base manifold
admits a K\"ahler-Einstein metric;
\item[] $f=K\t^{-1}+L$ for nonzero constants $K$ and $L$.
\end{list}
Conversely, in all even dimensions $n\geq 4$, there exist \sk\ pairs $(g,\t)$ and corresponding nonzero constants
$K$ and $L$ for which $(g/\t^2, K\t^{-1}+L)$ is a nontrivial quasi-Einstein pair.
\end{Theorem}
Note that the proof, relying on the local classification of \sk\ pairs, yields families, given explicitly, of all possible metrics $g$ (and hence $g/\t^2$) in Theorem \ref{A}, up to biholomorphic isometries.

The theory of \sk\ pairs has been used to obtain other results on K\"ahler metrics conformal to distinguished metrics. Namely, a similar result
holds if $g/\t^2$ is assumed to be Einstein (\cite{local}, see also \cite{pp} for examples), without, of course, any suppositions or conclusions on $f$.
Similarly, a completely analogous result holds if $g/\t^2$ is a gradient \sol\ \cite{confsol}. The present \bee\ case differs from the soliton case in two
regards: to carry out the classification in the latter case, it is enough to consider the soliton function $f=\t^{-1}$ (i.e. $K=1$, $L=0$); more importantly, in the soliton case, $g/\t^2$ is itself K\"ahler with respect to an oppositely oriented complex structure. Neither of these properties hold for the quasi-Einstein metrics in Theorem \ref{A}. On the other hand, in both the soliton and the \bee\ cases, the assumption that $\t$ gives rise to a Killing vector field can be dropped if one requires in advance that $f$ is affine in $\t^{-1}$.

We stress that while Theorem \ref{A} is of a local nature, the theory of \sk\ pairs also yields
a classification for compact manifolds. In fact, one has
\begin{Theorem}\lb{B}
With all assumptions on $M$, $g$, $\t$ and $f$ as in Theorem \ref{A}, assume also that $M$ is compact.
Then either $M$ is biholomorphic to complex projective space, or to the total space of a $CP^1$ bundle whose base manifold
admits a K\"ahler-Einstein metric.
\end{Theorem}
This follows directly from global \sk\ pair theory, specifically \cite[Theorem 29.2]{skrp}. Note further that
using local \sk\ theory, one obtains explicit families for the metric $g$. To obtain an exact classification of which of these metrics
in fact extends to a compact manifold, one must carry out further work analyzing boundary and positivity conditions on the function $Q=g(\nab\t,\nab\t)$,
along the lines of \cite{global}. It is reasonable to expect that the quasi-Einstein metrics appearing in such a global classification will
include those of L\"u, Page and Pope \cite{lpp} (because, for example, the compact spaces on which their metrics live are among those named in Theorem \ref{B}). If this is indeed the case, it will follow that their metrics are conformally K\"ahler.

The proof of Theorem \ref{A} proceeds as follows. Consideration of conformal changes yields that the pair $(g,\t)$ satisfies a Ricci-Hessian equation, with coefficients that force it to be a ``standard" \sk\ pair. The construction of such a pair depends on the existence of a nontrivial horizontal eigenfunction
for the Hessian of $\t$. This function is locally a function of $\t$, and is obtained by solving a pair of linear second
order odes. The desired nontrivial solution is obtained when $L/K$ is determined explicitly from constants
associated with the \sk\ pair, in such a way that it is nonzero. Inserting this solution as an ingredient in the canonical construction of an \sk\
pair yields the desired examples.

Note that Equations \Ref{rati} were obtained using a symbolic computation program.


\setcounter{Theorem}{0}
\renewcommand{\theTheorem}{\thesection.\arabic{Theorem}}

\section{\Bers\ and Ricci-Hessian pairs}
\setcounter{equation}{0}

On a manifold $M$ of real dimension larger than two, a pair $(g,\t)$ consisting of a
Riemannian metric $g$ and a smooth nonconstant function $\t$ is called a {\em Ricci-Hessian
pair} \cite{confsol} if the equation \be\lb{ric-hessian} \al\nab d\t +\ri=\gamma\, g\end{equation}
holds for the Hessian of $\t$, the Ricci tensor $\ri$ of $g$, and
some $C^\infty$ coefficient functions $\al$ and $\gamma$, possibly defined only on an open set
of $M$.

The pair $(g,f)$ satisfying \Ref{ber}, with $f$ nonconstant, is an example of a Ricci-Hessian pair.
Another example, and one that will be of main interest in what follows, is provided by the requirement that
a conformal change of $g$ yields a particular \ber. Namely,
suppose that $g$ is a Riemannian metric and $\t$, $f$ are nonconstant functions for which
the pair $(\wht{g}=g/\t^2,f)$ is a \ber. Thus $(-a/f)\wht{\nabla} df+\wht{r}=\lambda\wht{g}$ holds
for a constant $\lambda$ and a positive constant $a$, where $\wht{r}$, $\wht{\nabla}$
are the Ricci form and covariant derivative, respectively, of $\wht{g}$. Suppose further that
$df\we d\t=0$, so that $f$ is locally a function of $\t$. Using the formulas
$\wht{r} = r + (n - 2) \t^{-1}\nabla d\t + [\t^{-1}\Lap\t - (n - 1) \t^{-2}Q]g$ and
$\wht{\nabla} df = \nabla df + \t^{-1}[d\t\otimes df+ df\otimes d\t - g(\nabla\t,\nabla f)g]$
(see \cite[(2.1),(2.3)]{confsol}), a calculation as in \cite[(2.9)]{confsol} gives
\begin{eqnarray}\lb{f(tau)}
\ri&+&\left((n-2)\,\t^{-1}-a\,f'/f\right)\,\nab d\t\,
-(a/f)\left(f''+2\t^{-1}f'\right)\,d\t\otimes d\t\nonumber\\
&=&\left[\lambda\,\t^{-2}-
\t^{-1}\Delta\t+\left((n-1)\t^{-2}-af'\t^{-1}f^{-1}\right)Q\right]g.
\end{eqnarray}
Here $Q=g(\nab\t,\nab\t)$ and the prime denotes differentiation with respect to $\t$. If $f$ is affine in $\t^{-1}$,
the coefficient of $d\t\otimes d\t$ vanishes, and so we have:
\begin{Proposition}\lb{confqE}
Suppose $\t$ is a nonconstant smooth function and $g$ is a Riemannian metric that is conformal to a nontrivial \ber\  metric $g/\t^2$ for which the \berf\ is affine in the square root of the conformal factor (i.e. in $\t^{-1}$). Then $(g,\t)$ satisfies a Ricci-Hessian equation. The conclusion also follows if the \berf\ is merely locally a function of $\t$, provided that $g$ is K\"ahler and $\t$ is a Killing potential.
\end{Proposition}
Here $\t$ is a Killing potential if $J\nab\t$ is a Killing vector field, with $J$ the almost complex structure on $M$.
The last part of Proposition \ref{confqE} follows since if $f$ is locally a function of $\t$, the K\"ahler and Killing assumptions
imply that it is {\em in fact} affine in $\t^{-1}$. The proof of this claim is identical to that of
\cite[Proposition 3.1]{confsol}). Note that conversely, if $g$ is K\"ahler and $\t$ is a smooth nonconstant function such that $(g,\t)$ satisfies the Ricci-Hessian equation \Ref{ric-hessian}, then it follows that $\t$ is a Killing potential on the support of $\al$ (see \cite[beginning of \S3.1]{confsol}).

We note that if $K\t^{-1}+L$, with $K\neq 0$, is a \berf, so is $\t^{-1}+L/K$, for the same metric. This allows one to consider only \berfs\ of the form $f=\t^{-1}+k$ for a constant $k$,
as we do from here on.
From \Ref{f(tau)}, we then have that $(g,\t)$ is a Ricci-Hessian pair with
\be\lb{ag}\al=[n-2+a/(1+k\t)]/\t, \quad \gamma=\lambda\,\t^{-2}-\t^{-1}\Delta\t+[a/(1+k\t)+n-1]\t^{-2}Q.\end{equation}
Of course, this value of $\al$, but with $a=0$, is the one that occurs when
$g$ is conformally Einstein.



\section{Relation to \sk\ pairs}
\setcounter{equation}{0}

To define the notion of an \sk\ pair, denote by $\Mt$  the complement,
in a manifold $M$, of the critical set of a smooth function $\t$.  Recall that for a
Killing potential $\t$ on a K\"ahler manifold $M$, the set $\Mt$ is
open and dense in $M$.
\begin{Definition}\lb{skrp-def}\cite{local}
A nonconstant Killing potential $\t$ on a K\"ahler manifold $(M,J,g)$ is
called a {\em special K\"ahler-Ricci potential} if, on the set $\Mt$,
all non-zero tangent vectors orthogonal to $\nab\t$ and $J\nab\t$ are
eigenvectors of both $\nab d\t$ and $\ri$. The pair $(g,\t)$ is then called an {\em \sk\ pair}.
\end{Definition}

The utility of the concept of an \sk\ pair lies in the fact that it is a classifiable
geometric structure \cite{local, skrp,global}. We review part of the theory in Section \ref{de}.
Here we show:
\begin{Proposition}\lb{sol-skrp-2}
On a manifold of real dimension at least six, if $(g,\t)$ is a pair with $g$ K\"ahler
and $\t$ a nonconstant Killing potential, while the associated pair $(g/\t^2, f(\t))$ is a
nontrivial \berp, then $(g,\t)$ is an \sk\ pair. This also holds in real dimension four,
provided that $Q=g(\nabla\t,\nabla\t)$ and $\Delta\t$ are locally functions of $\t$.
\end{Proposition}

\begin{proof}
By Proposition \ref{confqE} $g$ satisfies a Ricci-Hessian equation \[\al\nab d\t +\ri=\gamma\, g.\]
The paragraphs past this proposition indicate that $f$ is in fact affine in $\t^{-1}$, and for the
purpose of determining $g$, it suffices to assume $f=\t^{-1}+k$, with $k$ a constant.

The coefficients $\al$, $\gamma$ of the Ricci-Hessian equation are then given by \Ref{ag}, from which one
sees that $\al$ is manifestly a function of $\t$, i.e. $d\al\we d\t=0$. It then follows that also
$d\gamma\we d\t=0$, if either $n\geq 6$, or else $n=4$ and both $Q=g(\nabla\t,\nabla\t)$ and
$\Delta\t$ are locally functions of $\t$ (see \cite[Proposition 3.3 and the paragraph before it]{confsol}).
Finally, it follows from \Ref{ag} that $\al d\al\neq 0$,
except on the sets where $\t=(2-n-a)/(\left(n-2\right)\,k)$ or
$\t=\left(-(n-2+a)\pm\sqrt{a(n+a-2)}\right)/{((n-2)\,k)}$.

These observations imply that the pair $(g,\t)$ is an \sk\ pair.
This follows from \cite[Proposition 3.5]{confsol} away
from the above mentioned degeneracy sets, and it follows in
the entire $\t$-noncritical set $\Mt$ by an argument analogous
to \cite[Corollary 3.7]{confsol}.
\end{proof}

\section{\sk\ pair theory and associated differential equations}\lb{de}
\setcounter{equation}{0}
By \cite[Definition $7.2$, Remark $7.3$]{local},
the \sk\ condition on $(g,\t)$ is equivalent to the existence,
on $\Mt$, of an orthogonal decomposition
$TM={\mathcal{V}}\oplus {\mathcal{H}}$, with
${\mathcal{V}}=\mathrm{span}(\nab\t, J\nab\t)$, along with
four smooth functions  $\phi$, $\psi$, $\beta$, $\mu$
which are pointwise eigenvalues for either $\nab d\t$ or $\ri$,
i.e., they satisfy
\begin{equation}\lb{skrp}
\begin{array}{rclrcl}
 \nab d\t|_{\mathcal{H}}\ &=&\phi\, g|_{\mathcal{H}}, \quad
&\nab d\t|_{\mathcal{V}}\ &=&\psi\, g|_{\mathcal{V}},\\
 \ri|_{\mathcal{H}}\ &=&\beta\, g|_{\mathcal{H}}, \quad
&\ri|_{\mathcal{V}}\ &=&\mu\, g|_{\mathcal{V}}.
\end{array}
\end{equation}
This decomposition is also $\ri$- and $\nab d\t$-orthogonal.

\begin{Remark}\lb{nontriv}
By \cite[Lemma 12.5]{local}, $\phi$ either vanishes identically on
$\Mt$, or never vanishes there. In the former case {\em only},
$g$ is reducible to a local product of K\"ahler metrics near any point
(see \cite[Corollary 13.2]{local} and \cite[Remark 16.4]{skrp}).
In the latter case, we call $g$ a {\em nontrivial} \sk\ metric.
\end{Remark}
\begin{Remark}\lb{c-const}
For a nontrivial \sk\ metric, consider $c=\t-Q/(2\phi)$,
with $Q=g(\nab \t,\nab \t)$, and $\kappa=\mathrm{sgn}(\phi)(\Lap\t+\lam\,Q/\phi)$,
regarded as functions $\Mt\ra \mathbb{R}$. By \cite[Lemma 12.5]{local},
$c$ is constant on $\Mt$, and will be called the {\em \sk\ constant}.
In any complex dimension $m\geq 2$, we will call a nontrivial \sk\ metric
{\em standard} if $\kappa$ is constant (and also use ``standard \sk\ pair"
as a designation for $(g,\t)$). According to \cite[\S27, using (10.1)
and Lemma 11.1]{local}, $\kappa$ is in fact constant if $m>2$, so that the
designation ``standard" involves an extra assumption as compared with ``nontrivial"
only when $m=2$. 
\end{Remark}

The following proposition summarizes a number of results given in \cite{confsol}.
\begin{Proposition}\lb{sta-sk}
For any \sk\ pair $(g,\t)$, the function $\phi$, i.e.
the horizontal eigenvalue of $\nabla d\t$, is locally a $C^\infty$ function of
$\t$ on $\Mt$. Furthermore, the pair satisfies a Ricci-Hessian equation on
the open set where $\nab d\t$ is not a multiple of $g$. If the pair is standard,
the coefficients $\al$ and $\gamma$ of this Ricci-Hessian equation are
locally functions of $\t$, and the function $\phi$ satisfies the ordinary differential
equation
\begin{equation}\lb{mek}
(\t-c)^2\phi''+\,(\t-c)[m-(\t-c)\al]\phi'-\,m\phi\,=\,-\mathrm{sgn}(\phi)\kappa/2.
\end{equation}
at points of $\Mt$ for which  $\phi'(\t)$ is nonzero. Also, on the same set,
\begin{equation}\lb{gamma}
\gamma=\al\phi+\left(\al(\t-c)-(m+1)\right)\phi'-(\t-c)\phi''
\end{equation} holds.
\end{Proposition}
\begin{proof}
The function $\phi$ is locally a function of $\t$ by \cite[Lemma 11.1a]{local}.
The existence of a Ricci-Hessian equation is the second part of
\cite[Proposition 3.5]{confsol}. Next, $\al$ is locally a function of $\t$ by
\cite[Remark 3.10]{confsol}, while the same then follows for $\gamma$ by
\cite[Proposition 3.3]{confsol}.  Equation \Ref{mek} holds by \cite[Proposition 4.1]{confsol},
and the expression \Ref{gamma} is obtained in \cite[second paragraph of Section 4.2]{confsol}.
\end{proof}
We will need the following refinement of this proposition, which extends the
domain of definition of the above equations, and was implicitly assumed in \cite{confsol}.
\begin{Proposition}\lb{refine}
Any Ricci-Hessian equation satisfied by an \sk\ pair $(g,\t)$ on some set, coincides
with the one of Proposition \ref{sta-sk} on the intersection of their domains.
Furthermore, Equations \Ref{mek} and \Ref{gamma} hold on the intersection of the (entire)
domain of such a Ricci-Hessian equation with $\Mt$, if the pair is nontrivial. These two equations
are an  ode, and, respectively, an expression for $\gamma$ in terms of functions of $\t$,
if the Ricci-Hessian equation has coefficients that are functions of $\t$, and the
pair is standard.
\end{Proposition}
\begin{proof}
The coefficients $\al$, $\gamma$ of the Ricci-Hessian equation in Proposition \ref{sta-sk}
are uniquely determined on the open set of $\Mt$ where $\nab d\t$ is not a multiple of $g$.
This follows as $\ri$ is a unique linear combination of $\nab d\t$  and $g$ on this set,
as the latter two tensors form at each point of this set a basis for the space of
twice covariant tensors for which the nonzero vectors ${\mathcal{H}}$ and ${\mathcal{V}}$ are eigenvectors
(see \cite[Remark 7.4]{local} and \cite[Remark 3.6]{confsol}). Thus any other Ricci-Hessian
equation must coincide with the above one, on the intersection of their domains.

Equations \Ref{mek} and \Ref{gamma} hold on any subset of $\Mt$ where a Ricci-Hessian equation holds,
as their derivation depends only on the relation $\beta-\mu=(\psi-\phi)\al$, which holds
at all points where the Ricci-Hessian equation holds
(see the proofs in \cite[Remark 3.10 and Proposition 4.1]{confsol}).
The final statement is immediate.
\end{proof}

We return now to our standard assumptions, as in, e.g.,  Proposition \ref{sol-skrp-2}:
$M$ is a manifold of real dimension at least six,  $(g,\t)$ is a pair with $g$ K\"ahler
and $\t$ a nonconstant Killing potential, while the associated pair $(g/\t^2, f(\t))$ is a
nontrivial \berp. $(g,\t)$ forms an \sk\ pair, by Proposition \ref{sol-skrp-2}.
{\em Assume also that $g$ is not reducible as a local product of K\"ahler metrics
in any neighborhood of a given point.} In other words, $g$ is nontrivial, and by
our assumption on the dimension, $(g,\t)$ is standard (see Remarks \ref{nontriv} and \ref{c-const}).
Proposition \ref{confqE} guarantees the existence of a Ricci-Hessian equation for the pair $(g,\t)$,
with coefficients given by \Ref{ag}, defined on the set where $\t\ne 0$. This set is open and dense
in $\Mt$, as $\t$ is a Killing potential, hence a Morse-Bott function. By Proposition \ref{refine},
Equations \Ref{mek} and \Ref{gamma} hold on the intersection of $\Mt$ and $\{\t\ne 0\}$ (both of which
are open dense sets in $M$).

To consider these equations explicitly, we substitute in Equation \Ref{mek}
the expression for $\al$  given in \Ref{ag}. This yields, after multiplying by $\t (1+k\t)$,
setting $n=2m$ and distributing various terms, the equation
\begin{eqnarray}\lb{1}
&&\t(\t-c)^2(1+k\t)\phi''+\,\left[(2-m)k\t^3+\left(2-m-a+(3m-4)kc\right)\t^2\right.\\
&+&\left.\left((2a+3m-4)c-2(m-1)kc^2\right)\t-(2m-2+a)c^2\right]\phi'
-\,\left[m\t+mk\t^2\right]\phi\,\nonumber\\
&=&\,-\mathrm{sgn}(\phi)\kappa\t (1+k\t)/2.\nonumber
\end{eqnarray}
Similarly, the expression \Ref{gamma} for $\gamma$ (again with the value of $\al$ from \Ref{ag}), must equal
its expression given in \Ref{ag} on $\{\t\ne 0\}$. Equating the two expressions, while
using $Q=2(\t-c)\phi$ (by definition of $c$) and $\Delta\t=2m\phi+2(\t-c)\phi'$ (\cite[(4.3.c)]{confsol})
gives
\begin{eqnarray}\lb{2}
&&\t^2(\t-c)(1+k\t)\phi''+\,\left[(1-m)k\t^3+(1-m-a+2mkc)\t^2\right.\\
&+&\left.c(a+2m)\t\right]\phi'
+\left[\left(a-2c(2m-1)k\right)\t-2c(a+2m-1)\right]\phi\,=\,-\lambda (1+k\t).\nonumber
\end{eqnarray}
We proceed to analyze these equations.

\section{Solutions of the equations}
\setcounter{equation}{0}
To examine the solutions of the system \Ref{1}--\Ref{2}, we add the product of \Ref{1} by $\t$
to the product of \Ref{2} by $-(\t-c)$. This results in a first order equation, which, after
rearrangement of terms yields, together with \Ref{1}, the system
\begin{eqnarray}\lb{qet}
&&\t(\t-c)(\t-2c)(1+k\t)\phi'+\left[-mk\t^3-\left(m+a-2c(2m-1)k\right)\t^2\right.\\
&+&\left.\left(c(3a+4m-2)-2c^2(2m-1)k\right)\t-2c^2(a+2m-1)\right]\phi\nonumber\\
&=&(1+k\t)\left[-\mathrm{sgn}(\phi)\kappa\t^2/2+\lambda(\t-c)\right],\nonumber\\
&&\t(\t-c)^2(1+k\t)\phi''+\,\left[(2-m)k\t^3+\left(2-m-a+(3m-4)kc\right)\t^2\right.\lb{qet2}\\
&+&\left.\left((2a+3m-4)c-2(m-1)kc^2\right)\t-(2m-2+a)c^2\right]\phi'
-\,\left[m\t+mk\t^2\right]\phi\,\nonumber\\
&=&\,-\mathrm{sgn}(\phi)\kappa\t(1+k\t)/2.\nonumber
\end{eqnarray}

To study solutions of this system, we recall the following (\cite[Lemma 4.3]{confsol}
\begin{Lemma}\lb{nonrat}
Let $\{\phi'+p\phi=q,\, A\phi''+B\phi'+C\phi=D\}$ be a system of ordinary differential
equations in the variable $\t$, with coefficients $p$, $q$, $A$, $B$, $C$ and $D$
that are rational functions. Then, on any nonempty interval admitting a solution
$\phi$, either
\be\lb{rat-rel} A(p^2-p')-Bp+C=0
\end{equation} holds identically, or
\be\lb{sol-doub}\phi=\left(D-A(q'-pq)-Bq\right)/\left(A(p^2-p')-Bp+C\right).
\end{equation} holds away from the (isolated) singularities of the right hand side.
\end{Lemma}

In applying Lemma \ref{nonrat} to the system formed by \Ref{qet}
and \Ref{qet2}, we of course modify \Ref{qet} appropriately, dividing it
by the factor $\t(\t-c)(\t-2c)(1+k\t)$. The resulting system has
a solution set identical to that of \Ref{1}--\Ref{2}
(certainly on intervals not containing $0$, $c$, $2c$ and $-1/k$ (if $k\neq 0$), and by a continuity
argument, on any interval). Computing \Ref{rat-rel} and the numerator of
\Ref{sol-doub} in this case using a symbolic computation program (and also verified by hand for the
case $k=0$), we get
\begin{equation}\lb{rati}
\begin{array}{lcl}\arraycolsep5pt
D-A(q'-pq)-Bq &=& 0,\\
A(p^2-p')-Bp+C&=&a(\t-c)^2(2ck+1)/((\t-2c)(\t k+1)).
\end{array}
\end{equation}
This immediately gives
\begin{Proposition}\lb{bc}
Suppose $a(2ck+1)\ne 0$. Then the system \Ref{1}--\Ref{2} has no nonzero solutions
on any nonempty open interval.
\end{Proposition}
\begin{proof}
Assume $a(2ck+1) \ne 0$. Then the right hand side of the second of Equations \Ref{rati} does not vanish
identically on the given interval, and thus so does the left hand side. Hence Lemma \ref{nonrat} implies
that any solution to the system \Ref{qet}--\Ref{qet2} is the ratio of the left-hand sides of the
two equations in \Ref{rati}, away from the point $c$. This ratio is the zero function. By continuity,
neither the system \Ref{qet}--\Ref{qet2}, nor the equivalent system formed by \Ref{1} and \Ref{2}, admits any
nonzero solutions on the given interval.
\end{proof}

\section{Solutions for the case $k=-1/(2c)$}
\setcounter{equation}{0}
If $k=-1/(2c)$, Equations \Ref{1}--\Ref{2} take, after multiplying by $2c$ and simplifying,
the form,
\begin{eqnarray}\lb{solsys}
&&\t(\t-c)^2(2c-\t)\phi''+\,\left[(m-2)\t^3+c\left(8-5m-2a\right)\t^2\right.\\
&+&\left.2c^2\left(2a+4m-5\right)\t-2c^3(2m-2+a)\right]\phi'
+\,\left[m\t(\t-2c)\right]\phi\,\nonumber\\
&=&\,\mathrm{sgn}(\phi)\kappa\t(\t-2c)/2,\nonumber\\
&&\t^2(\t-c)(2c-\t)\phi''+\,\left[(m-1)\t^3+2c(1-2m-a)\t^2\right.\nonumber\\
&+&\left.2c^2(a+2m)\t\right]\phi'
+\left[2c(a+2m-1)(\t-2c)\right]\phi\,=\,\lambda (\t-2c).\nonumber
\end{eqnarray}
We consider these solutions for $m$ a positive integer,  positive $a$, and {\em nonzero}
$c$ (as $2ck+1=0$, neither $c$ nor $k$ are zero).

One notices that these equations admit constant nonzero special solutions: $\kappa/(2m)$ with $\kappa>0$
for the first equation, and $\lambda/(2c(a+2m-1))$ for the second. Such a constant will be a
joint solution if these two values are equal (implying that $\lambda/c>0$).

A basis of joint solutions to the associated homogeneous equations is
obtained as follows. Each such solution must also solve the homogeneous first order
equation associated with \Ref{qet}, which is equivalent to $\phi'+p\phi=0$, where
in the expression for $p$ given before Equation \Ref{rati} one sets $k=-1/(2c)$, i.e.
\[p=\fr{a-1}{\t-2c}+ \fr m{\t-c}+\fr{1-a-2m}{\t}.\] The solutions of
this equation are constant multiples of $\exp(-\int p)$, i.e. of
\[(\t-2c)^{1-a}(\t-c)^{-m}\t^{2m-1+a}.\] To show that this solution indeed solves
the differential equations in \Ref{solsys}, say the first one, note that substituting
$\exp(-\int p)$ into $A\phi''+B\phi'+C\phi$ yields $\exp(-\int p)(A(p^2-p')-Bp+C)$, which
vanishes since \Ref{rat-rel} vanishes by the second equation in \Ref{rati} (with $2ck+1=0$).

To summarize, the general solution to the system \Ref{solsys} has the form
\be\lb{phi-sol}
\phi=C_1+C_2(\t-2c)^{1-a}(\t-c)^{-m}\t^{2m-1+a}
\end{equation}
for an arbitrary constant $C_2$, and a constant $C_1$ equal to both $\kappa/(2m)$ and
$\lambda/(2c(a+2m-1))$.


\section{Local geometry of standard \sk\ pairs}\lb{geo}
\setcounter{equation}{0}

As in \cite{confsol}, we review here the main case in the geometric classification of \sk\ metrics.
Let $\pi:(\hat{L},\langle\cdot,\cdot\rangle)\ra (N,h)$ be a Hermitian holomorphic line bundle
over a K\"ahler-Einstein manifold of complex dimension $m-1$. Assume that the
curvature of $\langle\cdot,\cdot\rangle$ is a multiple of the K\"ahler form of $h$.
Note that, if $N$ is compact and $h$ is not Ricci flat, this implies that $\hat{L}$ is smoothly
isomorphic to a rational power of the anti-canonical bundle of $N$.

Consider, on $\hat{L}\ssminus N$ (the total space of $\hat{L}$ excluding the zero section),
the metric $g$ given by
\be\lb{metric-form}
g|_{\mathcal{H}}=2|\t-c|\,\pi^*h, \quad g|_{\mathcal{V}}=\fr {Q(\t)}{(br)^2}\,
\mathrm{Re}\,\langle\cdot,\cdot\rangle,
\end{equation}
where\\
-- $\mathcal{V},\mathcal{H}$ are the vertical/horizontal distributions of $\hat{L}$, respectively,
the latter determined via the Chern connection of $\langle\cdot,\cdot\rangle$,\\
-- $c$ and $b\neq 0$ are constants,\\
-- $r$ is the norm induced by $\langle\cdot,\cdot\rangle$,\\
-- $\t=\t(r)$ is a function on $\hat{L}\ssminus N$, obtained by composing with $r$ another
function, denoted via abuse of notation by $\t$, and obtained as follows:
one fixes an open interval $I$ and a positive $C^\infty$ function $Q(\t)$ on $I$,
solves the differential equation $(b/Q)\,d\t=d(\log r)$ to obtain a diffeomorphism
$r(\t):I\ra (0,\infty)$, and defines $\t(r)$  as the inverse of this diffeomorphism.

The pair $(g,\t)$, with $\t=\t(r)$, is a nontrivial \sk\ pair (see \cite[\S8 and
\S16]{local}, as well as \cite[\S4]{skrp}), $|\nab\t|^2_g=Q(\t(r))$ holds, and
the connection on $\hat{L}$ is not flat.
The constant $\kappa$ of Remark \ref{c-const} is the Einstein constant of $h$,
so that as $g$ is nontrivial, it follows from the stipulation of a K\"ahler-Einstein base
that it is in fact standard (this is opposed to the case of an arbitrary \sk\
pair, in which $h$ need not be Einstein if $m=2$). For any $g$ standard, or merely nontrivial,
the \sk\ constant $c$ (see again Remark \ref{c-const}) coincides with $c$ of
\Ref{metric-form}.

Conversely, for any standard \sk\ pair $(g,\t)$ on a complex manifold $(M,J)$,
any point in the $\t$-noncritical set $\Mt$ has a neighborhood which is biholomorphically
isometric to an open set in some triple $(\hat{L}\ssminus N,g,\t(r))$ as above.
This claim is a special case of
\cite[Theorem 18.1]{local}). The biholomorphic isometry identifies
$\mathrm{span}\,(\nab\t, J\nab\t)$ and its orthogonal complement, with $\mathcal{V}$
and, respectively, $\mathcal{H}$. Moreover, whenever one can extend some $(g,\t(r))$ to all of $\hat{L}$,
such a biholomorphic isometry can also be defined on neighborhoods of points
in $M\ssminus \Mt$ \cite[Remark 16.4]{skrp}.

\section{Proof of Theorem \ref{A}}
\setcounter{equation}{0}
The proof of Theorem \ref{A} can now be concluded as follows. If the manifold $M$ has dimension
at least six, $g$ is K\"ahler and $(g/\t^2,f(\t))$ is a nontrivial \berp, then by Proposition \ref{sol-skrp-2},
$(g,\t)$ is an \sk\ pair. Since $g$ is K\"ahler-irreducible, by Remark \ref{nontriv} this \sk\ pair
is nontrivial, and with the assumption on the dimension, it is standard, according to Remark \ref{c-const}.
Thus, according to the classification given in Section \ref{geo}, $M$ is biholomorphic to an open set
in the total space of a $CP^1$ bundle whose base manifold admits a K\"ahler-Einstein metric
(and in fact much more is known about this bundle).

Following the paragraph past Proposition \ref{confqE}, the fact that $g$ is K\"ahler, $\t$ is Killing and
$f$ is locally a function of $\t$ implies that $f$ is affine in $\t^{-1}$, i.e. $f=K\t^{-1}+L$.
The constant $K$ is nonzero since $(g/\t^2,f)$ is a {\em nontrivial} \berp. Now replace
$f$ by $\t^{-1}+k:=\t^{-1}+L/K$, which is another \berf\ for the same metric $g/\t^2$.
According to \cite[Lemma 11.1a]{local}, the horizontal eigenfunction $\phi$ of $\nab d\t$
is locally a function of $\t$, and by Proposition \ref{refine} and the last paragraphs of Section \ref{de}, it satisfies the
system \Ref{1}--\Ref{2} at points of the $\t$-noncritical set $M_\t$  for which $\t$ is nonzero.
Using again Remark \ref{nontriv}, $\phi$ is nowhere vanishing on  $M_\t$. This together with Proposition \ref{bc} implies,
since $a\neq 0$, that $2ck+1=0$, for the constant $c$ associated to the \sk\ pair $(g,\t)$ as in Remark
\ref{c-const}. This implies $k\neq 0$, and hence $L$ is nonzero.

Finally, for any $m\geq 2$, fix a choice of data
$\pi:(\hat{L},\langle\cdot,\cdot\rangle)\ra (N,h)$, $c$ {\em nonzero} and $b$ satisfying all the criteria
required in defining a (standard) \sk\ pair, and also fix a constant $a>0$.
Using these values of $m$, $c$, and $a$, form the system of equations \Ref{solsys}, with $\lambda$
chosen so that the system admits constant solutions. Choose now some solution for $\phi$ of the form \Ref{phi-sol}.
For this solution, define $Q(\t)=2(\t-c)\phi(\t)$, and choose an interval $I$ where $Q$ is positive.
On this interval, solve $(b/Q)\,d\t=d(\log r)$ to obtain $\t(r)$ as the inverse of $r(\t)$.
Using $Q(\t(r))$, $\t(r)$ and the other data, obtain a standard \sk\ pair $(g,\t)$, with $\t$ given by $\t(r)$, and
the metric $g$ given by \Ref{metric-form}. Set $k=-1/(2c)$ and form $\al$ according to \Ref{ag},
and $\gamma$ by  expression \Ref{gamma} (with the value of $\al$ just described, and the chosen solution $\phi$).
As $\phi(\t)$ is a solution of \Ref{solsys},
it is also a solution for the equivalent system \Ref{1}--\Ref{2} with the choice of $k$ above. By the construction of
Equation \Ref{2}, the expression just obtained for $\gamma$ equals that in \Ref{ag} (with $Q$ and $\Lap\t$
given using their expressions just before \Ref{2}, which are valid for any nontrivial \sk\ pair). As \Ref{ag}
holds for the \sk\ pair just constructed, so does \Ref{f(tau)} for $f=(\t(r))^{-1}-1/(2c)$. This means exactly that
for the \sk\ pair $(g,\t(r))$, the pair $(g/(\t(r))^2, (\t(r))^{-1}-1/(2c))$ is a nontrivial quasi-Einstein
pair. This concludes the proof of Theorem \ref{A}.


\appendix
\section{Local obstructions to the existence of K\"ahler \bers}

In \cite{csw} it was shown that there exist no nontrivial K\"ahler \bers\
on compact manifolds. This was based on a structure theorem given for complete nontrivial
K\"ahler \bers\ on simply connected manifolds. While the latter theorem involves global
assumptions, the proof is mainly local. Here, using the same methods employed
in the proof of Theorem \ref{A},  we describe, under an extra hypothesis on
the dimension, an alternative approach to this result. As it is based on the theory of \sk\
pairs, it has the merit of leading to explicit expressions for the metric.

\begin{Theorem} On a manifold $M$ of complex dimension $m>2$, suppose $g$ is a K\"ahler
metric and $f$ a nonconstant function such that $(g,f)$ is a (nontrivial) \berp. Then
$g$ is reducible as a local product of K\"ahler metrics, one of whose components is
K\"ahler-Einstein, and the other situated on a two dimensional manifold, and given explicitly.
\end{Theorem}
\begin{proof}
The pair $(g,f)$ satisfies the Ricci-Hessian Equation \Ref{ric-hessian}, with
$\al=-a/f$ and $\gamma$ equal to the constant $\lambda$. Hence $d\al\we df=0$, $d\gamma\we df =0$
and $\al d\al\neq 0$ on points of the $f$-noncritical set $M_f$ where $f\neq 0$. Thus by
\cite[Proposition 3.5]{confsol}, $(g,f)$ is an \sk\ pair (even for $m=2$), and so $M$ is biholomorphic
to an open set in the total space of a $CP^1$ bundle whose base manifold admits a K\"ahler-Einstein metric.
We show that the \sk\ pair cannot be nontrivial. If it is nontrivial, as $m\geq 3$ it is standard, and so by Proposition
\ref{refine}, the ode \Ref{mek} (with the above  expression substituted for $\al$), holds for the
(nonzero) horizontal eigenfunction $\phi$ of $\nab d f$, at points
of $M_f$ for which $f$ is nonzero. Similarly, Equation \Ref{gamma} also holds there.
The resulting system is
\begin{eqnarray}\lb{Kqe}
f(f-c)^2\phi''&+&\,(f-c)[mf+(f-c)a]\phi'-\,mf\phi\,=\,-\mathrm{sgn}(\phi)\kappa f/2.\\
-f(f-c)\phi''&+&\left[-a(f-c)-(m+1)f\right]\phi'-a\phi=\lambda f.\nonumber
\end{eqnarray}
Adding the first equation to $f-c$ times the second gives the first order equation
\begin{equation}\lb{Kqe2}
-f(f-c)\phi'-(a(f-c)+mf)\phi=f(\lambda (f-c)-\mathrm{sgn}(\phi)\kappa/2).\end{equation}
Calculating the left hand side of \Ref{rati} for the pair consisting of the second equation
in \Ref{Kqe} and Equation \Ref{Kqe2} yields
$A(p^2-p'))-Bp+C=-a(f-c)/f$, $D-A(q'-pq)-Bq= 0$. Thus, as $a> 0$, there are no nowhere vanishing solutions $\phi$,
by Lemma \ref{nonrat}, and hence the \sk\ pair cannot be nontrivial.
If the \sk\ pair is trivial, it follows from the classification in \cite[Theorem 18.1]{local} of such pairs, that $g$ is reducible in the manner described in the theorem. Furthermore, the metric is given explicitly by a formula analogous to \Ref{metric-form}, namely $g|_{\mathcal{H}}=\pi^*h, \quad g|_{\mathcal{V}}=\fr {Q(\t)}{(br)^2}\,
\mathrm{Re}\,\langle\cdot,\cdot\rangle$, with $h$ the K\"ahler-Einstein metric and $Q$ having an explicit expression in $\t$ (see \cite[(19.1)]{local}.
\end{proof}
Note that the only possibility not explored above is that of a nontrivial \sk\ pair which is not standard, a
possibility which can occur only if $m=2$. One can investigate this further within \sk\ theory; however,
this possibility is ruled out by the result in \cite{csw}.



\end{document}